\documentclass[11pt]{article}

\usepackage{amsmath,amssymb,amsfonts,amsthm}
\usepackage{geometry}
\usepackage{graphicx}
\usepackage{placeins}
\usepackage[hidelinks]{hyperref}
\usepackage{mathrsfs}

\theoremstyle{remark}
\newtheorem{remark}{Remark}

\usepackage{bm}

\title{
Homogeneous Bracket Flows with Three-Form Flux on the Full Flag Manifold
\texorpdfstring{$\mathrm{SU}(3)/T^2$}{SU(3)/T2}
}

\author{
Francisco V. de Oliveira$^{1}$\thanks{francisco.vieira@ufersa.edu.br} \and
Francisco Brito$^{2}$\thanks{ fabrito@df.ufcg.edu.br}
}

\date{
$^{1}$ UFERSA-UFCG\\
$^{2}$ UFCG
}

\theoremstyle{plain}
\newtheorem{theorem}{Theorem}[section]
\newtheorem{proposition}[theorem]{Proposition}

\begin{document}

\hypersetup{
  pdftitle={Homogeneous Bracket Flows with Three-Form Flux on the Full Flag Manifold SU(3)/T2},
  pdfauthor={Francisco V. de Oliveira and Francisco Brito}
}
\maketitle

\begin{abstract}
We develop an algebraic framework for homogeneous metric flows with
closed three-form flux on reductive homogeneous spaces. The geometric
data are encoded by the complete reductive bracket
$\mu$ on $\mathfrak g=\mathfrak k\oplus\mathfrak m$, while the flux
is treated as an additional invariant datum $H$ satisfying
$d_\mu H=0$. For the metric equation
\[
\partial_tg
=
-2\operatorname{Ric}_g
+
\frac12\mathcal K_{g,H},
\]
we derive the corresponding bracket evolution
\[
\dot\mu
=
-\pi(P_{\mu,H})\mu,
\qquad
P_{\mu,H}
=
\operatorname{Ric}_{\mu}^{\sharp}
-
\frac14\mathcal K_{\mu,H}^{\sharp}.
\]

We apply this formulation to the full flag manifold
$\mathrm{SU}(3)/T^2$. For a prescribed closed invariant three-form,
we compute explicitly the quadratic flux operator, the induced
three-parameter metric system, and the evolution of the transported
bracket coefficients. The system has exactly four positive
stationary points. They correspond to the normal homogeneous
Einstein metric and the three invariant K\"ahler--Einstein metrics.
The normal stationary point has unstable index two, whereas each
K\"ahler--Einstein stationary point has unstable index one.

When the flux also evolves by the Hodge heat equation, we derive the
closed system for its amplitude and the three metric coefficients.
Every isotropic solution becomes extinct in finite time, with
$h^2/\lambda^2\to0$. For a nonzero transverse mode the absolute
perturbation is of order $(T-t)^{3/5}$, whereas the perturbation
relative to the collapsing scale is of order $(T-t)^{-2/5}$.
We relate the invariant three-form to the NS--NS sector of the strict
$\mathrm{SU}(3)$-structure ansatz in type IIA supergravity.
\end{abstract}



\section{Introduction}

Generalized flag manifolds provide homogeneous internal spaces on
which invariant metrics, curvature and differential forms can be
described by Lie algebra data \cite{Caviezel2008,GasparimGramaVarea2023}.
Scherk--Schwarz reductions illustrate a related algebraic viewpoint:
a constant duality twist determines the gauging in the reduced theory
\cite{ScherkSchwarz1979,Hull2002,Samtleben2008}. On a flag manifold,
root brackets encode the internal homogeneous geometry. The comparison
concerns the role of algebraic data in the two constructions.

Fusi, Lafuente and Stanfield \cite{FusiLafuenteStanfield2024}
developed a bracket-flow formulation of generalized Ricci flow on
discrete quotients of Lie groups. Here we work with the complete
reductive bracket of a homogeneous space with isotropy. On the full
flag manifold $\mathrm{SU}(3)/T^2$, this formulation retains the
$\mathfrak{k}$-component of $[\mathfrak{m},\mathfrak{m}]$ and
allows us to compare prescribed flux with flux evolving by the
Hodge heat equation.

We formulate homogeneous metric flows with a closed invariant
three-form using the complete reductive bracket, including its
isotropy components. On the full flag manifold $\mathrm{SU}(3)/T^2$,
we compute the Ricci endomorphism, construct an exact invariant
three-form $\eta_0$ and determine its quadratic contribution to the
metric equation. We then study the same metric system in two regimes.
For prescribed $H=h\eta_0$ with $h\ne0$, the flow selects a unique
scale on each of the four invariant Einstein rays. We classify the
stationary points and compute their linear stability. When $H$ evolves
by the Hodge heat equation, we derive the equation for $h(t)$ and solve
the isotropic system. All isotropic solutions become extinct in finite
time. Their linearized transverse perturbations grow relative to the
isotropic scale even as the absolute perturbations decay. The
Einstein directions themselves are classical
\cite{PediconiSbiti2022}; our prescribed-flux calculation fixes the
scale on each ray and gives the spectrum of the unnormalized metric
system at the selected metric. The three-form used here is defined
by an invariant primitive, whereas \cite{GasparimGramaVarea2023}
construct a three-form by infinitesimal T-duality.

Section~\ref{sec:homogeneous-bracket-flow} gives the general bracket
formulation. Section~\ref{sec:su3-bracket-flow} derives the explicit
system shared by the prescribed and evolving flux regimes, studied
in Sections~\ref{sec:prescribed-flux} and~\ref{sec:evolving-flux}.
Section~\ref{sec:type-IIA-compactification} places the invariant data
in a type IIA ansatz and discusses the requirements for interpreting
the flow in the lower-dimensional scalar theory. The geometric-flow
approach to infinite-distance limits motivates that discussion
\cite{BaumePalti2016,CalderonInfante2020,LustPaltiVafa2019,
Kehagias2020,DemulderRaml2025}.

\section{Homogeneous metric flows and reductive brackets}
\label{sec:homogeneous-bracket-flow}

We use Lauret's homogeneous bracket formulation
\cite{Lauret2013CurvatureFlowsHomogeneous,Lauret2011BracketFlow}.
For a bracket-flow treatment of generalized Ricci flow on discrete
quotients of Lie groups, see \cite{FusiLafuenteStanfield2024}.

The metric and three-form equations below are related to geometric
flows with torsion \cite{GarciaFernandez2013,StreetsTian2010,GarciaFernandezStreets2021}.
\subsection{Homogeneous brackets and invariant fluxes}
\label{subsec:homogeneous-brackets-fluxes}

Let
\[
\mathfrak g
=
\mathfrak k\oplus\mathfrak m
\]
be a reductive decomposition associated with a homogeneous space
$G/K$. The homogeneous structure is encoded by the complete bracket
\[
\mu
\in
\Lambda^2\mathfrak g^*\otimes\mathfrak g,
\]
satisfying the Jacobi identity and the reductive conditions
\[
\mu(\mathfrak k,\mathfrak k)\subseteq\mathfrak k,
\qquad
\mu(\mathfrak k,\mathfrak m)\subseteq\mathfrak m.
\]

In general, one has
\[
\mu(\mathfrak m,\mathfrak m)
\subseteq
\mathfrak k\oplus\mathfrak m.
\]
Consequently, the projection
\[
\mu_{\mathfrak m}
\in
\Lambda^2\mathfrak m^*\otimes\mathfrak m
\]
does not, by itself, encode the complete homogeneous geometry. In particular, the isotropy projection is nonzero in the full flag
example below.

A homogeneous flux background is described by a pair
\[
(\mu,H),
\]
where
\[
H\in\Lambda^3\mathfrak m^*
\]
is $\operatorname{Ad}(K)$-invariant and closed with respect to the
differential induced by the homogeneous bracket:
\[
d_\mu H=0.
\]

Thus, a closed invariant three-form is an additional geometric datum.
It is not determined by the bracket in complete generality. In the
flowing-flags correspondence, one may construct distinguished fluxes
from infinitesimal T-duality and from the root data of the associated
nilradicals. Such a construction defines a particular assignment
\[
\mu\longmapsto H(\mu),
\]
but this assignment must be specified and its compatibility with the
evolution must be verified.

For root-compatible invariant forms, the nonzero components of $H$
may be supported on admissible root triples. Schematically, one may
write
\[
H
=
\sum_{\alpha,\beta}
h_{\alpha,\beta}\,
e^\alpha\wedge e^\beta\wedge e^{-(\alpha+\beta)},
\]
where the sum is taken over the triples for which
\[
\alpha+\beta\in R.
\]
The coefficients $h_{\alpha,\beta}$ depend on the flux construction
and should not, without further justification, be identified with the
root constants $N_{\alpha,\beta}$.

\subsection{Transport of the metric and flux}

Fix a reference inner product $Q$ on $\mathfrak m$ and let $\mu_0$
be the initial homogeneous bracket in a fixed presentation of $G/K$.
Write the metric equation as
\begin{equation}
\partial_t g=-2g(P_{g,H}\,\cdot,\cdot),
\qquad
P_{g,H}=\operatorname{Ric}^{\sharp}_g-\tfrac14\mathcal K^\sharp_{g,H}.
\label{eq:generalized-ricci-flow}
\end{equation}
Here $\mathcal K_{g,H}(X,Y)=\sum_{p,q}H(X,e_p,e_q)H(Y,e_p,e_q)$
for a $g$-orthonormal basis $\{e_p\}$, and the sharp denotes raising
one index with $g$. Thus \eqref{eq:generalized-ricci-flow} is
$\partial_t g=-2\operatorname{Ric}_g+\tfrac12\mathcal K_{g,H}$.

Choose a $K$-equivariant isometry
$A(t)\colon(\mathfrak m,g(t))\to(\mathfrak m,Q)$ solving
$\dot A=-A P_{g,H}$. Such a choice follows by starting from a
$K$-equivariant isometry at the initial time and differentiating
$g=A^*Q$. Extend $A$ by the identity on $\mathfrak k$ and set
\begin{equation}
\mu=A\cdot\mu_0,\qquad \widehat H=(A^{-1})^*H,\qquad
(A\cdot\mu_0)(X,Y)=A\mu_0(A^{-1}X,A^{-1}Y).
\label{eq:transported-homogeneous-bracket}
\end{equation}
Then $d_\mu\widehat H=(A^{-1})^*d_{\mu_0}H$, and $P_{\mu,\widehat H}
=A P_{g,H}A^{-1}$ on $\mathfrak m$, extended by zero on
$\mathfrak k$. Differentiating the action yields
\begin{equation}
\dot\mu=-\pi(P_{\mu,\widehat H})\mu,\qquad
(\pi(L)\mu)(X,Y)=L\mu(X,Y)-\mu(LX,Y)-\mu(X,LY).
\label{eq:bracket-flow-with-flux}
\end{equation}
In components in a fixed basis, the equation reads
\[
\dot\mu_{BC}^{D}
=-P^D{}_E\mu_{BC}^{E}+P^E{}_B\mu_{EC}^{D}+P^E{}_C\mu_{BE}^{D}.
\]

The transported form satisfies
\begin{equation}
\partial_t\widehat H
=(A^{-1})^*(\partial_t H)-\pi_3(P_{\mu,\widehat H})\widehat H,
\qquad
(\pi_3(L)\omega)(X_1,X_2,X_3)
=-\sum_{i=1}^3\omega(X_1,\ldots,LX_i,\ldots,X_3).
\label{eq:transported-flux}
\end{equation}
Thus a fixed form on $G/K$ can have time-dependent coefficients in
the transported frame. If the flux evolves in the fixed presentation,
its evolution law supplies the first term of
\eqref{eq:transported-flux}. If one instead prescribes an algebraic
assignment $\widehat H=\widehat H(\mu)$, its derivative along
\eqref{eq:bracket-flow-with-flux} must agree with
\eqref{eq:transported-flux} for the chosen flux law.

\section{Homogeneous bracket flow on \texorpdfstring{$\mathrm{SU}(3)/T^2$}{SU(3)/T2}}
\label{sec:su3-bracket-flow}

\subsection{Reductive decomposition and invariant brackets}
\label{subsec:su3-reductive-bracket}

Let
\[
\mathfrak g=\mathfrak{su}(3),
\qquad
\mathfrak k=\mathfrak t^2,
\]
and consider the $Q$-orthogonal reductive decomposition
\[
\mathfrak g
=
\mathfrak k\oplus\mathfrak m,
\qquad
\mathfrak m
=
\mathfrak m_1\oplus\mathfrak m_2\oplus\mathfrak m_3,
\]
where $Q=-B_{\mathfrak{su}(3)}$, with $B(X,Y)=6\operatorname{tr}(XY)$, and
\[
\mathfrak m_1=\mathfrak m_{\alpha_1},
\qquad
\mathfrak m_2=\mathfrak m_{\alpha_2},
\qquad
\mathfrak m_3=\mathfrak m_{\alpha_3}.
\]

Each $\mathfrak m_i$ is a real two-dimensional irreducible
$\operatorname{Ad}(T^2)$-module. The bracket relations are
\[
[\mathfrak k,\mathfrak m_i]\subseteq\mathfrak m_i,
\qquad
[\mathfrak m_i,\mathfrak m_i]\subseteq\mathfrak k,
\]
and, for pairwise distinct indices $\{i,j,k\}=\{1,2,3\}$,
\[
[\mathfrak m_i,\mathfrak m_j]\subseteq\mathfrak m_k.
\]

Here $\alpha_1=e_1-e_2$, $\alpha_2=e_2-e_3$ and $\alpha_3=\alpha_1+\alpha_2$. With $E_{ab}$
denoting the elementary matrices, set
\[
A_{ab}=\frac{E_{ab}-E_{ba}}{\sqrt{12}},\qquad
S_{ab}=\frac{\mathrm i(E_{ab}+E_{ba})}{\sqrt{12}} .
\]
We take $(X_i,Y_i)=(A_{ab},S_{ab})$ for $(a,b)=(1,2),(2,3),(1,3)$,
respectively. These are $Q$-orthonormal bases of the three summands.

\subsection{Metric-dependent bracket coefficients}
\label{subsec:metric-dependent-bracket}

An $\mathrm{SU}(3)$-invariant Riemannian metric on
$\mathrm{SU}(3)/T^2$ is determined by three positive parameters:
\[
g_{\boldsymbol\lambda}
=
\lambda_1 Q|_{\mathfrak m_1}
+
\lambda_2 Q|_{\mathfrak m_2}
+
\lambda_3 Q|_{\mathfrak m_3},
\qquad
\lambda_i>0.
\]

Define the positive symmetric endomorphism
\[
A_{\boldsymbol\lambda}|_{\mathfrak k}
=
I_{\mathfrak k},
\qquad
A_{\boldsymbol\lambda}|_{\mathfrak m_i}
=
\sqrt{\lambda_i}\,I_{\mathfrak m_i}.
\]
Then
\[
g_{\boldsymbol\lambda}(X,Y)
=
Q(A_{\boldsymbol\lambda}X,A_{\boldsymbol\lambda}Y),
\qquad
X,Y\in\mathfrak m.
\]

The corresponding bracket on the fixed inner-product space
$(\mathfrak g,Q)$ is
\[
\mu_{\boldsymbol\lambda}
=
A_{\boldsymbol\lambda}\cdot\mu_0,
\]
where
\[
(A\cdot\mu_0)(X,Y)=A\,\mu_0(A^{-1}X,A^{-1}Y).
\]

For the chosen vectors $X_i$, write
\[
[X_j,X_k]_{\mathfrak m_i}
=
c_iX_i,
\qquad
\{i,j,k\}=\{1,2,3\},
\]
where $c_i$ is determined by the chosen normalization of the
Cartan--Weyl basis. Then
\[
\mu_{\boldsymbol\lambda}(X_j,X_k)_{\mathfrak m_i}
=
c_i
\sqrt{\frac{\lambda_i}{\lambda_j\lambda_k}}\,X_i.
\]

We therefore introduce the metric-dependent bracket coefficients
\[
x_1
=
c_1\sqrt{\frac{\lambda_1}{\lambda_2\lambda_3}},
\qquad
x_2
=
c_2\sqrt{\frac{\lambda_2}{\lambda_1\lambda_3}},
\qquad
x_3
=
c_3\sqrt{\frac{\lambda_3}{\lambda_1\lambda_2}}.
\]
For the standard symmetric normalization of the $A_2$ root system,
one may choose
\[
|c_1|=|c_2|=|c_3|=:c,
\]
with the signs determined by the orientation of the root vectors.

The isotropy components scale differently. If
\[
[X_i,Y_i]_{\mathfrak k}=K_i
\]
then
\[
\mu_{\boldsymbol\lambda}(X_i,Y_i)_{\mathfrak k}
=
\frac{1}{\lambda_i}K_i.
\]
Consequently, the complete transported bracket contains both
\[
x_i
\sim
\sqrt{\frac{\lambda_i}{\lambda_j\lambda_k}}
\]
and isotropy coefficients proportional to $\lambda_i^{-1}$.

This shows that the homogeneous bracket is not described by a
single evolving parameter. Its diagonal invariant sector is
parametrized by the three metric variables
\[
(\lambda_1,\lambda_2,\lambda_3),
\]
or, equivalently, by the associated bracket coefficients together
with the isotropy components.

\subsection{Ricci endomorphism}
\label{subsec:flag-ricci}

By isotropy equivariance, $\operatorname{Ric}^{\sharp}|_{\mathfrak m_i}
=r_i I_{\mathfrak m_i}$. For $Q=-B$, the squared bracket constant
is $[123]=1/3$ and $\dim\mathfrak m_i=2$. The homogeneous Ricci formula
gives
\begin{equation}
 r_i(\lambda)=\frac{1}{2\lambda_i}
 +\frac{1}{12}\left(
 \frac{\lambda_i}{\lambda_j\lambda_k}
 -\frac{\lambda_j}{\lambda_i\lambda_k}
 -\frac{\lambda_k}{\lambda_i\lambda_j}\right).
 \label{eq:flag-ricci-eigenvalues-corrected}
\end{equation}

The formula holds for $\{i,j,k\}=\{1,2,3\}$, with
$\operatorname{Ric}_{g_\lambda}|_{\mathfrak m_i}=r_i g_\lambda|_{\mathfrak m_i}$.
Appendix~\ref{app:ricci-su3} gives the normalization details.

\subsection{An invariant closed three-form}
\label{subsec:invariant-three-form-su3}

We now fix a nonzero closed
$\mathrm{SU}(3)$-invariant three-form on
\[
\frac{\mathrm{SU}(3)}{T^2}.
\]
This flux is regarded as an additional homogeneous datum and is not
assumed to coincide with the flux produced on the T-dual flag in the
flowing-flags correspondence.

For the chosen bases $\{X_i,Y_i\}$, let
\[
\{\theta_i,\varphi_i\}
\]
denote the corresponding dual coframe. Define the invariant
two-forms
\[
\omega_i
=
\theta_i\wedge\varphi_i.
\]

The differential of an invariant two-form
\[
B_0
=
b_1\omega_1+b_2\omega_2+b_3\omega_3
\]
is an invariant three-form.

The zero-weight relation $\alpha_1+\alpha_2-\alpha_3=0$ gives the
invariant complex product
\[
(\theta_1+i\varphi_1)\wedge(\theta_2+i\varphi_2)\wedge
(\theta_3-i\varphi_3).
\]
Its negative imaginary part is the real form
\begin{equation}
\begin{aligned}
\eta_0={}&
\theta_1\wedge\theta_2\wedge\varphi_3
-\theta_1\wedge\varphi_2\wedge\theta_3\\
&-\varphi_1\wedge\theta_2\wedge\theta_3
-\varphi_1\wedge\varphi_2\wedge\varphi_3.
\end{aligned}
\label{eq:invariant-eta-su3}
\end{equation}

With the chosen bases, the invariant differential satisfies
\begin{equation}
d\omega_1=d\omega_2=-d\omega_3=c_0\eta_0,
\qquad c_0=\frac1{\sqrt{12}}.
\label{eq:closed-flux-normalization}
\end{equation}
Indeed, the projected brackets from $\mathfrak m_1\times\mathfrak m_2$
are
\[
[X_1,X_2]=c_0X_3,\qquad [X_1,Y_2]=[Y_1,X_2]=c_0Y_3,\qquad
[Y_1,Y_2]=-c_0X_3.
\]
The other cyclic brackets
follow from $Q([U,V],W)=Q(U,[V,W])$. Inserting these brackets into the
invariant exterior differential gives \eqref{eq:closed-flux-normalization}.
In particular, $\eta_0=c_0^{-1}d\omega_1$ is exact.
Consequently,
\begin{equation}
dB_0
=
c_0(b_1+b_2-b_3)\eta_0.
\label{eq:dB0-su3}
\end{equation}

We use the ansatz
\begin{equation}
H=h\eta_0 .
\label{eq:H-ansatz-su3}
\end{equation}
The amplitude $h$ is constant in Section~\ref{sec:prescribed-flux}
and evolves in Section~\ref{sec:evolving-flux}. Exactness is preserved
in both cases. The quadratic flux tensor depends on the form itself.

\subsection{The quadratic flux operator}
\label{subsec:quadratic-flux-operator-su3}

Let
\[
g_{\boldsymbol\lambda}
=
\lambda_1Q|_{\mathfrak m_1}
+
\lambda_2Q|_{\mathfrak m_2}
+
\lambda_3Q|_{\mathfrak m_3}.
\]
A $g_{\boldsymbol\lambda}$-orthonormal basis of $\mathfrak m_i$ is
given by
\[
E_i=\frac{X_i}{\sqrt{\lambda_i}},
\qquad
F_i=\frac{Y_i}{\sqrt{\lambda_i}}.
\]
The corresponding orthonormal coframe is
\[
e^i=\sqrt{\lambda_i}\,\theta_i,
\qquad
f^i=\sqrt{\lambda_i}\,\varphi_i.
\]

Consequently,
\[
\theta_i=\frac{e^i}{\sqrt{\lambda_i}},
\qquad
\varphi_i=\frac{f^i}{\sqrt{\lambda_i}},
\]
and the flux becomes
\begin{equation}
\begin{aligned}
H=\frac{h}{\sqrt{\lambda_1\lambda_2\lambda_3}}
\Bigl(
&e^1\wedge e^2\wedge f^3
-e^1\wedge f^2\wedge e^3\\
&-f^1\wedge e^2\wedge e^3
-f^1\wedge f^2\wedge f^3
\Bigr).
\end{aligned}
\label{eq:H-orthonormal-coframe}
\end{equation}

Set
\[
\sigma(\boldsymbol\lambda,h)
=
\frac{h}{\sqrt{\lambda_1\lambda_2\lambda_3}}.
\]
Each orthonormal vector $E_i$ or $F_i$ occurs in exactly two
orthonormal components of $H$. Since the definition
\[
\mathcal K_{\boldsymbol\lambda,H}(X,Y)
=
\sum_{p,q}
H(X,e_p,e_q)H(Y,e_p,e_q)
\]
sums over ordered pairs $(p,q)$, each component contributes twice.
Therefore,
\[
\mathcal K_{\boldsymbol\lambda,H}(E_i,E_i)
=
\mathcal K_{\boldsymbol\lambda,H}(F_i,F_i)
=
4\sigma(\boldsymbol\lambda,h)^2.
\]
All mixed terms vanish:
\[
\mathcal K_{\boldsymbol\lambda,H}(E_i,F_i)=0,
\qquad
\mathcal K_{\boldsymbol\lambda,H}
(\mathfrak m_i,\mathfrak m_j)=0
\quad
(i\neq j).
\]

It follows that
\begin{equation}
\mathcal K_{\boldsymbol\lambda,H}^{\sharp}
|_{\mathfrak m_i}
=
\frac{4h^2}
{\lambda_1\lambda_2\lambda_3}
I_{\mathfrak m_i},
\qquad
i=1,2,3.
\label{eq:K-operator-su3}
\end{equation}
Hence,
\begin{equation}
\kappa_1
=
\kappa_2
=
\kappa_3
=
\frac{4h^2}
{\lambda_1\lambda_2\lambda_3}.
\label{eq:kappa-su3}
\end{equation}

The modified Ricci eigenvalues are therefore
\begin{equation}
q_i(\boldsymbol\lambda,h)
=
r_i(\boldsymbol\lambda)
-
\frac{h^2}
{\lambda_1\lambda_2\lambda_3}.
\label{eq:modified-ricci-eigenvalues-su3}
\end{equation}

\subsection{The common metric system}
\label{subsec:flag-metric-system}

Write $F_i(\boldsymbol\lambda,h)$ for the right-hand side of
\begin{equation}
\dot\lambda_i=F_i(\boldsymbol\lambda,h)
=-2r_i(\boldsymbol\lambda)\lambda_i+\frac{2h^2}{\lambda_j\lambda_k},
\qquad \{i,j,k\}=\{1,2,3\}.
\label{eq:flux-metric-general}
\end{equation}
Inserting \eqref{eq:flag-ricci-eigenvalues-corrected} gives
\begin{equation}
\begin{aligned}
\dot\lambda_1&=-1-\frac16\left(
\frac{\lambda_1^2}{\lambda_2\lambda_3}
-\frac{\lambda_2}{\lambda_3}-\frac{\lambda_3}{\lambda_2}\right)
+\frac{2h^2}{\lambda_2\lambda_3},\\
\dot\lambda_2&=-1-\frac16\left(
\frac{\lambda_2^2}{\lambda_1\lambda_3}
-\frac{\lambda_1}{\lambda_3}-\frac{\lambda_3}{\lambda_1}\right)
+\frac{2h^2}{\lambda_1\lambda_3},\\
\dot\lambda_3&=-1-\frac16\left(
\frac{\lambda_3^2}{\lambda_1\lambda_2}
-\frac{\lambda_1}{\lambda_2}-\frac{\lambda_2}{\lambda_1}\right)
+\frac{2h^2}{\lambda_1\lambda_2}.
\end{aligned}
\label{eq:corrected-flux-system}
\end{equation}
When $h=0$, this reduces to the ordinary homogeneous Ricci flow.

Permutation symmetry makes the diagonal locus invariant for either
flux law. Along $\lambda_i=\lambda$ we have $r_i=5/(12\lambda)$ and
\begin{equation}
\dot\lambda=-\frac56+\frac{2h^2}{\lambda^2}.
\label{eq:isotropic-flux-flow}
\end{equation}

\subsection{Explicit evolution of the bracket coefficients}
\label{subsec:explicit-bracket-coefficients}

Since $\dot\lambda_i/\lambda_i=-2q_i$, logarithmic differentiation
of $x_i=c_i\sqrt{\lambda_i/(\lambda_j\lambda_k)}$ gives
$\dot x_i=-x_i(q_i-q_j-q_k)$. Substituting
\eqref{eq:modified-ricci-eigenvalues-su3} gives
\begin{equation}
\dot x_i
=
-x_i
\left(
r_i-r_j-r_k
+
\frac{h^2}{\lambda_1\lambda_2\lambda_3}
\right).
\label{eq:bracket-coefficient-flow-general}
\end{equation}

For $i=1$, the Ricci combination is
\[
r_1-r_2-r_3
=
\frac{1}{2\lambda_1}
-
\frac{1}{2\lambda_2}
-
\frac{1}{2\lambda_3}
+
\frac{1}{12}
\left(
\frac{3\lambda_1}{\lambda_2\lambda_3}
-
\frac{\lambda_2}{\lambda_1\lambda_3}
-
\frac{\lambda_3}{\lambda_1\lambda_2}
\right).
\]
Thus,
\begin{equation}
\begin{aligned}
\dot x_1
=
-x_1\Biggl[
&
\frac{1}{2\lambda_1}
-
\frac{1}{2\lambda_2}
-
\frac{1}{2\lambda_3}\\
&+
\frac{1}{12}
\left(
\frac{3\lambda_1}{\lambda_2\lambda_3}
-
\frac{\lambda_2}{\lambda_1\lambda_3}
-
\frac{\lambda_3}{\lambda_1\lambda_2}
\right)
+
\frac{h^2}{\lambda_1\lambda_2\lambda_3}
\Biggr].
\end{aligned}
\label{eq:x1-explicit-flow}
\end{equation}
The equations for $x_2$ and $x_3$ are obtained by cyclic
permutation of the indices.

The isotropy components of the transported bracket have coefficients
proportional to
\[
y_i=\frac{1}{\lambda_i}.
\]
Their evolution is
\begin{equation}
\dot y_i
=
2q_i y_i
=
2y_i
\left(
r_i
-
\frac{h^2}{\lambda_1\lambda_2\lambda_3}
\right).
\label{eq:isotropy-coefficient-flow}
\end{equation}

Consequently, the complete reductive bracket evolution contains both
the cross-summand coefficients $x_i$ and the isotropy coefficients
$y_i$. They evolve simultaneously according to the same modified
Ricci eigenvalues $q_i$.

In the transported frame the three-form is
$\widehat H=h(\lambda_1\lambda_2\lambda_3)^{-1/2}\eta_0$, where the
coframe defining $\eta_0$ is now the fixed reference coframe. The
coefficient formulas above apply to both choices of $h$.

\section{Stationary brackets for prescribed flux}
\label{sec:prescribed-flux}

Throughout this section $H=h\eta_0$ is fixed on the homogeneous space
and $h\ne0$. We classify the equilibria of
\eqref{eq:corrected-flux-system} and then compare their linearizations.
\subsection{Classification of the positive stationary points}

We now determine all stationary points of the prescribed-flux
system in the positive cone
\[
\mathcal P
=
\left\{
(\lambda_1,\lambda_2,\lambda_3)\in\mathbb R^3:
\lambda_i>0
\right\}.
\]

Let
\[
F(\lambda_1,\lambda_2,\lambda_3)
=
(F_1,F_2,F_3)
\]
be the vector field defined by
\eqref{eq:corrected-flux-system}. A direct computation gives
\begin{equation}
F_1-F_2
=
\frac{
(\lambda_1-\lambda_2)
\left[
\lambda_3^2-(\lambda_1+\lambda_2)^2+12h^2
\right]
}{
6\lambda_1\lambda_2\lambda_3
}.
\label{eq:difference-F1-F2}
\end{equation}
Cyclically,
\begin{align}
F_2-F_3
&=
\frac{
(\lambda_2-\lambda_3)
\left[
\lambda_1^2-(\lambda_2+\lambda_3)^2+12h^2
\right]
}{
6\lambda_1\lambda_2\lambda_3
},
\\
F_3-F_1
&=
\frac{
(\lambda_3-\lambda_1)
\left[
\lambda_2^2-(\lambda_3+\lambda_1)^2+12h^2
\right]
}{
6\lambda_1\lambda_2\lambda_3
}.
\end{align}

We first show that every positive stationary point has at least two
equal components. Suppose, by contradiction, that
\[
\lambda_1,\lambda_2,\lambda_3
\]
are pairwise distinct. Since
\[
F_1=F_2=F_3=0,
\]
the three difference identities imply
\begin{align}
\lambda_3^2
&=
(\lambda_1+\lambda_2)^2-12h^2,
\\
\lambda_1^2
&=
(\lambda_2+\lambda_3)^2-12h^2,
\\
\lambda_2^2
&=
(\lambda_3+\lambda_1)^2-12h^2.
\end{align}
Subtracting the second identity from the first gives
\[
\lambda_3^2-\lambda_1^2
=
(\lambda_1+\lambda_2)^2
-
(\lambda_2+\lambda_3)^2.
\]
Factoring both sides, we obtain
\[
(\lambda_3-\lambda_1)(\lambda_3+\lambda_1)
=
-(\lambda_3-\lambda_1)
(\lambda_1+\lambda_3+2\lambda_2).
\]
Since $\lambda_3\neq\lambda_1$, division by
$\lambda_3-\lambda_1$ gives
\[
\lambda_3+\lambda_1
=
-(\lambda_1+\lambda_3+2\lambda_2),
\]
which is impossible because all $\lambda_i$ are positive. Therefore,
at least two components must be equal.

Suppose, without loss of generality,
that
\[
\lambda_1=\lambda_2=a,
\qquad
\lambda_3=b.
\]
The stationary equations reduce to
\[
b^2-6ab+12h^2=0
\]
and
\[
4a^2+b^2=12h^2.
\]
Eliminating $h^2$ gives
\[
(b-a)(b-2a)=0.
\]

The first possibility, $b=a$, yields the isotropic stationary point
\[
a=b=\sqrt{\frac{12}{5}}\,|h|.
\]
The second possibility, $b=2a$, gives
\[
a^2=\frac32h^2,
\qquad
b^2=6h^2.
\]
Therefore,
\[
a=\sqrt{\frac32}\,|h|,
\qquad
b=2\sqrt{\frac32}\,|h|.
\]

We have consequently proved the following result.

\begin{proposition}
Let $h\neq0$. The prescribed-flux metric flow on
$\mathrm{SU}(3)/T^2$ has exactly four stationary points in the
positive cone. One is isotropic:
\[
\boldsymbol\lambda_{\mathrm{iso}}
=
\sqrt{\frac{12}{5}}\,|h|\,(1,1,1).
\]
The remaining three are anisotropic:
\[
\begin{aligned}
\boldsymbol\lambda_{\mathrm{aniso}}^{(1)}
&=
\sqrt{\frac32}\,|h|\,(2,1,1),\\
\boldsymbol\lambda_{\mathrm{aniso}}^{(2)}
&=
\sqrt{\frac32}\,|h|\,(1,2,1),\\
\boldsymbol\lambda_{\mathrm{aniso}}^{(3)}
&=
\sqrt{\frac32}\,|h|\,(1,1,2).
\end{aligned}
\]
Thus, the anisotropic stationary metrics have ratios given by the
three permutations of $(1,1,2)$.
\end{proposition}

\begin{remark}
At a stationary point, the equation
\[
q_i(\boldsymbol\lambda,h)=0
\]
implies
\[
r_1=r_2=r_3
=
\frac{h^2}{\lambda_1\lambda_2\lambda_3}.
\]
Consequently,
\[
\operatorname{Ric}_{\boldsymbol\lambda}^{\sharp}
=
\frac{h^2}{\lambda_1\lambda_2\lambda_3}
I_{\mathfrak m},
\]
and every stationary metric of the prescribed-flux flow is an
invariant Einstein metric.

The four stationary points found above correspond precisely to the
four invariant Einstein metrics on $\mathrm{SU}(3)/T^2$: the normal
homogeneous metric with ratio
\[
(1,1,1)
\]
and the three K\"ahler--Einstein metrics with ratios
\[
(1,1,2),\qquad
(1,2,1),\qquad
(2,1,1).
\]
Thus, the prescribed flux does not create new invariant Einstein
directions. Instead, it selects a definite scale along each of the
four classical Einstein rays.
\end{remark}

\subsection{Linear stability}

\paragraph{Isotropic point.}
Let
\[
F(\lambda_1,\lambda_2,\lambda_3)
=
(F_1,F_2,F_3)
\]
denote the vector field defined by
\eqref{eq:corrected-flux-system}. At the isotropic stationary
point
\[
\boldsymbol\lambda_*
=
(\lambda_*,\lambda_*,\lambda_*),
\qquad
\lambda_*^2=\frac{12}{5}h^2,
\]
the Jacobian matrix is
\begin{equation}
DF(\boldsymbol\lambda_*)
=
-\frac{1}{3\lambda_*}
\begin{pmatrix}
1&2&2\\
2&1&2\\
2&2&1
\end{pmatrix}.
\label{eq:jacobian-isotropic-fixed-point}
\end{equation}

The radial direction
\[
v_{\mathrm{rad}}=(1,1,1)
\]
is an eigenvector with eigenvalue
\[
\nu_{\mathrm{rad}}
=
-\frac{5}{3\lambda_*}.
\]
Hence, the stationary point is attracting along the isotropic locus.

The anisotropic subspace is
\[
\mathcal A
=
\left\{
(v_1,v_2,v_3)\in\mathbb R^3:
v_1+v_2+v_3=0
\right\}.
\]
Every vector in $\mathcal A$ is an eigenvector with eigenvalue
\[
\nu_{\mathrm{aniso}}
=
\frac{1}{3\lambda_*},
\]
with multiplicity two. Therefore, anisotropic perturbations grow to
first order.

We conclude that the isotropic stationary metric is hyperbolic and
has unstable index two:
\begin{equation}
\dim E^s=1,
\qquad
\dim E^u=2.
\label{eq:stable-unstable-dimensions}
\end{equation}
In particular, the stationary point is asymptotically stable within
the isotropic invariant locus but is a saddle point in the full
three-parameter space of invariant metrics.

\paragraph{Anisotropic points.}

Consider
\[
\boldsymbol\lambda_{\mathrm{aniso}}
=
(a,a,2a),
\qquad
a=\sqrt{\frac32}\,|h|.
\]
At this stationary point, the Jacobian is
\begin{equation}
DF(a,a,2a)
=
-\frac{1}{6a}
\begin{pmatrix}
1&5&1\\
5&1&1\\
4&4&4
\end{pmatrix}.
\label{eq:jacobian-anisotropic-fixed-point}
\end{equation}
Its eigenvalues are
\begin{equation}
-\frac{4}{3a},
\qquad
-\frac{1}{3a},
\qquad
\frac{2}{3a}.
\label{eq:eigenvalues-anisotropic-fixed-point}
\end{equation}
Hence, the anisotropic stationary point is hyperbolic, with
\[
\dim E^s=2,
\qquad
\dim E^u=1.
\]
The other two anisotropic stationary points are obtained by
permuting the coordinates and have the same spectrum. Therefore,
all three anisotropic stationary points are saddles of unstable
index one.

\subsection{Stationarity of the complete reductive bracket}

At every metric equilibrium $q_i=0$ for all $i$. The cross-summand
and isotropy coefficients in
\eqref{eq:bracket-coefficient-flow-general} and
\eqref{eq:isotropy-coefficient-flow} are therefore stationary.
Conversely, stationarity of the complete bracket forces
$\dot y_i=0$ for all $i$, because each
$[X_i,Y_i]_{\mathfrak k}$ is nonzero. Since $y_i=\lambda_i^{-1}$,
this gives a metric equilibrium. The transported form is stationary
there as well.

\begin{figure}[htbp]
\centering
\includegraphics[width=.82\linewidth]{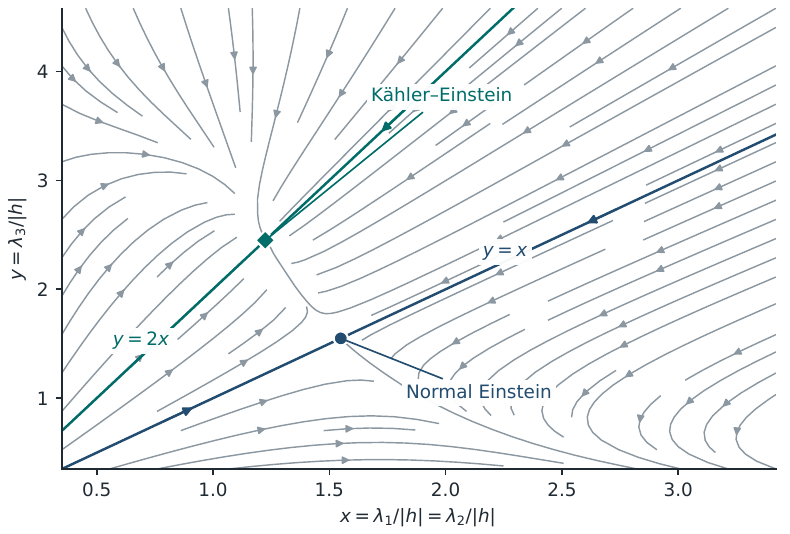}
\caption{Phase portrait of the prescribed-flux system on the invariant slice
$\lambda_1=\lambda_2$, in coordinates $x=\lambda_1/|h|=\lambda_2/|h|$ and
$y=\lambda_3/|h|$. The marked points represent the normal Einstein metric
and one of the three K\"ahler--Einstein metrics.}
\label{fig:prescribed-flux}
\end{figure}
\FloatBarrier

\section{Evolution of the three-form and isotropic collapse}
\label{sec:evolving-flux}

\subsection{The coupled invariant system}

On the fixed homogeneous space we couple \eqref{eq:flux-metric-general}
to
\begin{equation}
\partial_t H=-\Delta_{d,g}H,\qquad
\Delta_{d,g}=dd^*_g+d^*_gd.
\label{eq:flux-heat-law}
\end{equation}
Here $\Delta_{d,g}$ is nonnegative. Since $H$ is closed, the flux
equation is $\partial_t H=-dd^*_g H$. The codifferential of $\eta_0$
is invariant, so it lies in the span of the three orthogonal
invariant two-forms $\omega_i$. Their pointwise squared norms are
$\lambda_i^{-2}$, and the four orthogonal monomials in
\eqref{eq:invariant-eta-su3} give
$|\eta_0|^2_{g_\lambda}=4/(\lambda_1\lambda_2\lambda_3)$
for the exterior-power norm. By compactness and invariance, the
$L^2$ adjoint identity reduces to the pointwise equality
$\langle d^*_{g_\lambda}\eta_0,\omega_i\rangle
=\langle\eta_0,d\omega_i\rangle$. Consequently,
\[
d^*_{g_\lambda}\eta_0
=\frac{4c_0}{\lambda_1\lambda_2\lambda_3}
\left(\lambda_1^2\omega_1+\lambda_2^2\omega_2-\lambda_3^2\omega_3\right).
\]
Applying $d$ and using \eqref{eq:closed-flux-normalization} shows
that the ansatz $H=h(t)\eta_0$ is preserved. The amplitude satisfies
\begin{equation}
\dot h
=-\frac{1}{3}
\frac{\lambda_1^2+\lambda_2^2+\lambda_3^2}
{\lambda_1\lambda_2\lambda_3}\,h.
\label{eq:transverse-flux-system}
\end{equation}
Together with \eqref{eq:corrected-flux-system}, this defines a smooth
ODE on $(0,\infty)^3\times\mathbb R$. In the transported description
\eqref{eq:transported-flux}, the heat term is
$-\Delta_{d_\mu,Q}\widehat H$.

\subsection{Explicit isotropic solution and extinction}
\begin{theorem}
\label{thm:flag-isotropic-asymptotics}
Consider the isotropic generalized
Ricci flow with evolving exact
three-form $H=h(t)\eta_0$ and
$Q=-B_{\mathfrak{su}(3)}$.

Its evolution equations are
\begin{equation}
\begin{aligned}
\dot\lambda&=-\frac56+\frac{2h^2}{\lambda^2},\\
\dot h&=-\frac h\lambda.
\end{aligned}
\label{eq:normalized-isotropic-coupled}
\end{equation}

Every solution with $\lambda_0>0$
becomes extinct in finite forward time
$T$. Moreover,
\[
\frac{h(t)^2}{\lambda(t)^2}
\longrightarrow0,
\qquad
\lambda(t)\sim\frac56(T-t),
\qquad
h(t)\longrightarrow0.
\]
\end{theorem}

\begin{proof}
On the diagonal, \eqref{eq:transverse-flux-system} reduces to
\[
\dot h=-\frac h\lambda.
\]

Set
\[
z=\frac{h^2}{\lambda^2},
\qquad
\frac{d\tau}{dt}=\frac1\lambda.
\]
Then
\begin{equation}
\begin{aligned}
z'&=-\frac13z-4z^2,\\
(\log\lambda)'&=-\frac56+2z,\\
h'&=-h.
\end{aligned}
\label{eq:normalized-rescaled-system}
\end{equation}

Writing $z_0=h_0^2/\lambda_0^2$,
direct integration gives
\[
z(\tau)=
\frac{z_0e^{-\tau/3}}
{1+12z_0(1-e^{-\tau/3})},
\]
\[
\lambda(\tau)
=\lambda_0e^{-5\tau/6}
\sqrt{1+12z_0(1-e^{-\tau/3})},
\]
and
\[
h(\tau)=h_0e^{-\tau}.
\]

It follows that
\[
z(\tau)\longrightarrow0
\]
and
\[
\lambda(\tau)\sim
\lambda_0\sqrt{1+12z_0}\,
e^{-5\tau/6}.
\]
Thus
\[
T-t_0=\int_0^\infty
\lambda(\tau)\,d\tau<\infty.
\]

Finally,
\[
\dot\lambda=-\frac56+2z
\longrightarrow-\frac56,
\]
and hence
\[
\lambda(t)\sim\frac56(T-t).
\]
The explicit expression for $h$
also shows that $h(t)\to0$.
\end{proof}

\begin{remark}
The explicit solution gives $|h(t)|\asymp(T-t)^{6/5}$ if
$h_0\ne0$. Thus the flux term $2h^2/\lambda^2$ in the equation for
$\dot\lambda$ is of order $(T-t)^{2/5}$. With the exterior-power norm,
\[
|H|_{g_\lambda}^2=\frac{4h^2}{\lambda^3}
\asymp (T-t)^{-3/5}.
\]
The pointwise squared norm therefore diverges for nonzero initial flux,
while its ratio to $\operatorname{Scal}_{g_\lambda}=5/(2\lambda)$
tends to zero.
\end{remark}

\subsection{Transverse linear instability of
isotropic collapsing solutions}
\label{subsec:transverse-instability}

We linearize \eqref{eq:corrected-flux-system} and
\eqref{eq:transverse-flux-system} along the isotropic solutions.
Let $(\lambda(t),h(t))$ be an isotropic
solution, so that
\[
\lambda_1=\lambda_2=\lambda_3=\lambda.
\]
Consider infinitesimal anisotropic
perturbations
\[
\lambda_i=\lambda+\varepsilon_i,
\qquad
\varepsilon_1+\varepsilon_2+\varepsilon_3=0.
\]

\begin{proposition}
\label{prop:transverse-variational-equation}
The transverse variational equation along
an isotropic solution is
\begin{equation}
\dot\varepsilon_i
=
\left(
-\frac1{2\lambda}
+\frac{2h^2}{\lambda^3}
\right)\varepsilon_i.
\label{eq:transverse-variational}
\end{equation}
\end{proposition}

\begin{proof}
We use the difference identity \eqref{eq:difference-F1-F2}.
At the isotropic locus, the coefficient
of $\lambda_1-\lambda_2$ becomes
\[
-\frac1{2\lambda}
+\frac{2h^2}{\lambda^3}.
\]
The corresponding identities hold
cyclically.

By permutation symmetry, the anisotropic
subspace
\[
\mathcal A=
\left\{
(\varepsilon_1,\varepsilon_2,\varepsilon_3):
\sum_i\varepsilon_i=0
\right\}
\]
is invariant under the linearization.
The linearized flux equation does not
couple to $\mathcal A$, since its
derivatives with respect to the three
metric variables coincide at isotropy.
This proves the assertion.
\end{proof}

\begin{theorem}
\label{thm:relative-transverse-instability}
Let $(\lambda(t),h(t))$ be any isotropic
solution with $\lambda(t)>0$, and define
\[
\frac{d\tau}{dt}=\frac1{\lambda(t)},
\qquad \tau(t_0)=0.
\]
Every solution of the transverse
variational equation satisfies
\begin{equation}
\frac{\varepsilon_i(\tau)}{\lambda(\tau)}
=
\frac{\varepsilon_i(0)}{\lambda_0}
e^{\tau/3}.
\label{eq:relative-anisotropic-growth}
\end{equation}

In particular, every nonzero transverse
perturbation grows relative to the
isotropic metric scale as
$\tau\to\infty$.

Moreover, if $T$ denotes the finite
extinction time of the isotropic
background, then
\[
|\varepsilon_i(t)|
\asymp(T-t)^{3/5}
\]
for every nonzero initial component,
whereas
\[
\frac{\|\varepsilon(t)\|}{\lambda(t)}
\asymp(T-t)^{-2/5}.
\]
Thus the isotropic collapsing trajectory
is linearly unstable with respect to
relative anisotropic perturbations.
\end{theorem}

\begin{proof}
Set
\[
z=\frac{h^2}{\lambda^2}.
\]
In the rescaled time $\tau$, the
transverse equation becomes
\[
\varepsilon_i'
=\left(-\frac12+2z\right)\varepsilon_i.
\]
On the other hand, the isotropic metric
satisfies
\[
(\log\lambda)'=-\frac56+2z.
\]
Consequently,
\[
\left(\frac{\varepsilon_i}{\lambda}\right)'
=
\frac13\frac{\varepsilon_i}{\lambda},
\]
which integrates to
\eqref{eq:relative-anisotropic-growth}.

For the normalization $Q=-B$, the
isotropic solution is explicitly
\[
\lambda(\tau)=
\lambda_0e^{-5\tau/6}
\sqrt{1+12z_0(1-e^{-\tau/3})},
\qquad
z_0=\frac{h_0^2}{\lambda_0^2}.
\]
It follows that
\[
\varepsilon_i(\tau)=
\varepsilon_i(0)e^{-\tau/2}
\sqrt{1+12z_0(1-e^{-\tau/3})}.
\]
Thus
\[
\lambda(\tau)\asymp e^{-5\tau/6},
\qquad
|\varepsilon_i(\tau)|
\asymp e^{-\tau/2}
\]
for nonzero initial components.

Since
\[
\lambda(t)\sim\frac56(T-t),
\]
we have
\[
T-t\asymp e^{-5\tau/6}.
\]
Substitution yields
\[
|\varepsilon_i(t)|
\asymp(T-t)^{3/5}
\]
and
\[
\frac{\|\varepsilon(t)\|}{\lambda(t)}
\asymp(T-t)^{-2/5}.
\]
\end{proof}

\begin{remark}
The instability established above is a
property of the linearized evolution
relative to the collapsing background.
It does not, by itself, imply nonlinear
instability in a specified function
space or classify the singularities
of nearby anisotropic solutions.
Such conclusions require a separate
analysis of the full nonlinear system.
\end{remark}

\begin{figure}[htbp]
\centering
\includegraphics[width=\linewidth]{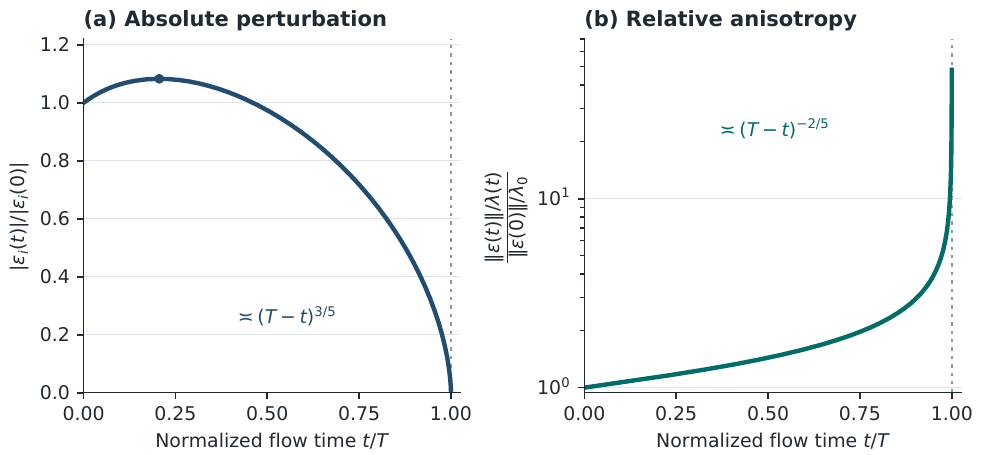}
\caption{Transverse linearized perturbations for $h_0^2/\lambda_0^2=1/2$.
Time is measured from the initial instant. The absolute perturbation
decays near extinction at order $(T-t)^{3/5}$, and the relative
anisotropy grows at order $(T-t)^{-2/5}$.}
\label{fig:evolving-flux}
\end{figure}
\FloatBarrier

\section{Relation to type IIA compactification on \texorpdfstring{$\mathrm{SU}(3)/T^2$}{SU(3)/T2}}
\label{sec:type-IIA-compactification}

We consider massive type IIA supergravity on a ten-dimensional
space-time of the form
\begin{equation}
M_{10}
=
M_4
\times
\frac{\mathrm{SU}(3)}{T^2},
\label{eq:ten-dimensional-product}
\end{equation}
where $M_4$ denotes the four-dimensional external space-time and
$\mathrm{SU}(3)/T^2$ is the six-dimensional internal flag manifold.

We place the invariant data in the NS--NS and R--R ansatz for massive
type IIA. The internal metric parameters can be promoted to
four-dimensional scalar fields. The geometric time of the preceding
sections is independent of the external spacetime coordinates.

\subsection{Ten-dimensional bosonic sector}
\label{subsec:type-IIA-bosonic-sector}

In the string frame, the NS--NS sector of type IIA supergravity is
described by
\begin{equation}
S_{\mathrm{NS}}
=
\frac{1}{2\kappa_{10}^2}
\int_{M_{10}}
e^{-2\Phi}
\left(
R_{10}\star_{10}1
+
4\,d\Phi\wedge\star_{10}d\Phi
-
\frac12H\wedge\star_{10}H
\right).
\label{eq:NSNS-action-string-frame}
\end{equation}
Equivalently,
\[
S_{\mathrm{NS}}
=
\frac{1}{2\kappa_{10}^2}
\int_{M_{10}}
d^{10}x\,
\sqrt{-g_{10}}\,
e^{-2\Phi}
\left(
R_{10}
+
4|d\Phi|^2
-
\frac{1}{12}H_{MNP}H^{MNP}
\right).
\]

In the democratic formulation, the R--R field strengths are
collected into
\[
F
=
F_0+F_2+F_4+F_6+F_8+F_{10},
\]
subject to the duality conditions
\[
F_p
=
(-1)^{\frac{p(p-1)}{2}}
\star_{10}F_{10-p}.
\]
The Romans mass is
\[
F_0=m.
\]

The total pseudo-action may be written schematically as
\begin{equation}
S_{\mathrm{IIA}}
=
S_{\mathrm{NS}}
-
\frac{1}{8\kappa_{10}^2}
\int_{M_{10}}
\sum_{p=0,2,\ldots,10}
F_p\wedge\star_{10}F_p
+
S_{\mathrm{loc}},
\label{eq:type-IIA-pseudo-action}
\end{equation}
where $\star_{10}$ is the string-frame Hodge operator and
$S_{\mathrm{loc}}$ denotes possible source terms.
The duality conditions supplement the democratic pseudo-action
\cite[Appendix~A]{Caviezel2008}.

\subsection{Homogeneous metric ansatz}
\label{subsec:IIA-metric-ansatz}

We adopt the ten-dimensional Einstein-frame metric ansatz
\begin{equation}
ds_{10,E}^2
=
e^{2A(x)}ds_{4,E}^2
+
e^{2B(x)}
ds_6^2(\boldsymbol\lambda),
\label{eq:ten-dimensional-metric-ansatz}
\end{equation}
where
\begin{equation}
ds_6^2(\boldsymbol\lambda)
=
\lambda_1Q|_{\mathfrak m_1}
+
\lambda_2Q|_{\mathfrak m_2}
+
\lambda_3Q|_{\mathfrak m_3}.
\label{eq:internal-flag-metric}
\end{equation}

The functions
\[
\lambda_i=\lambda_i(x),
\qquad
\Phi=\Phi(x),
\]
are four-dimensional scalar fields. The warp factors $A$ and $B$
will eventually be chosen so that $ds_{4,E}^2$ is in the
four-dimensional Einstein frame.

For the vacuum ansatz below, $A$, $B$, $\Phi$ and the $\lambda_i$
are constant on spacetime. Since $g_{10,s}=e^{\Phi/2}g_{10,E}$, the
physical string-frame internal coefficients are
$e^{2B+\Phi/2}\lambda_i$ in this parametrization. In the following
subsections we absorb this constant factor into the $\lambda_i$ and
use those symbols for the string-frame coefficients. The external
volume form in the vacuum flux formulas is also in the string frame.

\subsection{Invariant \texorpdfstring{$\mathrm{SU}(3)$}{SU(3)} structure}
\label{subsec:invariant-SU3-structure}

Use the coframe of Section~\ref{subsec:invariant-three-form-su3} and set
\begin{equation}
J=\lambda_1\omega_1+\lambda_2\omega_2-\lambda_3\omega_3,
\qquad
z^i=\sqrt{\lambda_i}(\theta_i+i\epsilon_i\varphi_i),
\qquad
(\epsilon_1,\epsilon_2,\epsilon_3)=(1,1,-1).
\label{eq:IIA-fundamental-form}
\end{equation}
For the strict $\mathrm{SU}(3)$ conventions of \cite{Caviezel2008},
take the complex volume form
\begin{equation}
\Omega=-i\,\overline{z^1\wedge z^2\wedge z^3}.
\label{eq:IIA-complex-volume-form}
\end{equation}
The zero-weight relation ensures invariance and the chosen phase gives
\begin{equation}
J\wedge\Omega=0,\qquad
-\frac{i}{8}\Omega\wedge\overline\Omega
=\frac{1}{3!}J^3=\operatorname{vol}_6,\qquad
\operatorname{Re}\Omega=\sqrt{\lambda_1\lambda_2\lambda_3}\,\eta_0.
\label{eq:SU3-structure-normalization}
\end{equation}

\subsection{Invariant flux ansatz}
\label{subsec:IIA-flux-ansatz}

The prescribed NS--NS flux considered in the geometric-flow sector
is
\begin{equation}
H
=
h\,\eta_0,
\qquad
d\eta_0=0.
\label{eq:IIA-H-flux}
\end{equation}
With the conventions introduced above, $\eta_0$ is proportional to
one of the real components of the invariant complex volume form.

For the R--R sector, the invariant vacuum ansatz is
\begin{align}
F_0
&=
m,
\\
F_2
&=
f_1\omega_1+f_2\omega_2+f_3\omega_3,
\\
F_4
&=
f_{\mathrm{ext}}\operatorname{vol}_4
+
e_1\widetilde\omega^{\,1}
+
e_2\widetilde\omega^{\,2}
+
e_3\widetilde\omega^{\,3},
\\
F_6^{\mathrm{int}}
&=
e_0\operatorname{vol}_6,
\label{eq:IIA-RR-flux-ansatz}
\end{align}
where
\[
\omega_i=\theta_i\wedge\varphi_i
\]
and $\widetilde\omega^{\,i}$ are invariant four-forms dual to
$\omega_i$, normalized by
\begin{equation}
\int_{\mathrm{SU}(3)/T^2}
\omega_i\wedge\widetilde\omega^{\,j}
=
\delta_i{}^j.
\label{eq:invariant-form-pairing}
\end{equation}

The invariant ansatz uses the coefficients
\[
m,\quad h,\quad f_i,\quad e_i,\quad e_0.
\]
Bianchi identities and any tadpole conditions constrain these data.
Flux quantization constrains the corresponding charges.
The coefficient $h$ of the exact NS--NS form does not label a
nonzero NS--NS cohomology class.

\subsection{Bianchi identities}
\label{subsec:IIA-Bianchi-identities}

In the absence of NS5-brane sources, the NS--NS flux satisfies
\begin{equation}
dH=0.
\label{eq:Bianchi-H}
\end{equation}
The R--R Bianchi identities are
\begin{equation}
dF_p+H\wedge F_{p-2}=j_{p+1},
\label{eq:IIA-Bianchi-general}
\end{equation}
where $j_{p+1}$ denotes the signed current of the corresponding
localized or smeared sources. In the notation of \cite{Caviezel2008},
$j_3^{\mathrm{D6/O6}}=-j_6$.

In particular,
\begin{align}
dF_0&=0,\\
dF_2+F_0H&=j_3,\\
dF_4+H\wedge F_2&=j_5,\\
dF_6+H\wedge F_4&=j_7.
\label{eq:IIA-Bianchi-components}
\end{align}

For D6/O6 sources, the relevant identity is
\begin{equation}
dF_2+mH=j_3^{\mathrm{D6/O6}}.
\label{eq:D6-O6-Bianchi}
\end{equation}
A source-free truncation requires
\[
dF_2=-mH.
\]
By \eqref{eq:closed-flux-normalization}, this is equivalent to
$c_0(f_1+f_2-f_3)=-mh$. The identity can always be solved at the
level of invariant real forms. Supersymmetry and flux quantization
impose further conditions.

In the democratic formulation, the external Freund--Rubin component
of $F_4$ and the internal component of $F_6$ are related by the
ten-dimensional duality constraint. Consequently,
$f_{\mathrm{ext}}$ and $e_0$ should not be regarded as independent
flux parameters. The full $F_6$ also contains a component of the form
$\operatorname{vol}_4\wedge\beta_2$, with $\beta_2$ fixed by duality
from the internal part of $F_4$. The identity for $dF_6$ uses the
full ten-dimensional form.

In the analysis of the ten-dimensional equations, we take
\[
F_0,\qquad F_2,\qquad F_4
\]
as the independent R--R field strengths. The higher-degree forms are
then reconstructed from the democratic duality relations.

\subsection{Strict \texorpdfstring{$\mathrm{SU}(3)$}{SU(3)} structure for \texorpdfstring{$\mathcal N=1$ AdS$_4$}{N=1 AdS4}}
\label{subsec:supersymmetric-IIA-ansatz}

We now specialize the invariant flux ansatz to the strict
$\mathrm{SU}(3)$-structure class of supersymmetric massive type IIA
AdS$_4$ compactifications. We assume a constant warp factor and a
constant dilaton $\Phi$.

In the conventions adopted here, the supersymmetric fluxes take the
form
\begin{align}
H
&=
\frac{2m}{5}e^{\Phi}\operatorname{Re}\Omega,
\label{eq:susy-H-flux}
\\
F_0
&=
m,
\label{eq:susy-F0}
\\
F_2
&=
\frac{f}{9}J+F_2',
\label{eq:susy-F2}
\\
F_4
&=
f\,\operatorname{vol}_4
+
\frac{3m}{10}J\wedge J,
\label{eq:susy-F4}
\end{align}
where $m$ is the Romans mass, $f$ is the Freund--Rubin parameter and
$F_2'$ is the primitive $(1,1)$ component of $F_2$:
\begin{equation}
J\lrcorner F_2'=0,
\qquad
F_2'\wedge J\wedge J=0.
\label{eq:primitive-F2}
\end{equation}

The torsion classes of the internal $\mathrm{SU}(3)$-structure obey
\begin{align}
dJ
&=
-\frac32 iW_1^-\operatorname{Re}\Omega,
\label{eq:dJ-torsion}
\\
d\Omega
&=
W_1^-J\wedge J
+
W_2^-\wedge J,
\label{eq:dOmega-torsion}
\end{align}
with
\begin{equation}
W_1^-
=
-\frac{4i}{9}e^\Phi f,
\qquad
W_2^-
=
-i e^\Phi F_2'.
\label{eq:torsion-flux-relations}
\end{equation}

For example, \eqref{eq:closed-flux-normalization} gives
$dJ=c_0(\lambda_1+\lambda_2+\lambda_3)\eta_0$. The first torsion
equation therefore fixes the Freund--Rubin parameter through
\[
e^\Phi f=-\frac{3c_0(\lambda_1+\lambda_2+\lambda_3)}
{2\sqrt{\lambda_1\lambda_2\lambda_3}}.
\]
The R--R Bianchi identity involving the primitive flux is
\[
dF_2+mH=j_3^{\mathrm{D6/O6}}.
\]
Using the supersymmetric flux ansatz, it reduces to an equation for
the primitive component:
\begin{equation}
dF_2'
=
\left(
\frac{2}{27}f^2
-
\frac{2}{5}m^2
\right)
e^\Phi\operatorname{Re}\Omega
+
j_3^{\mathrm{D6/O6}},
\label{eq:susy-Bianchi-F2prime}
\end{equation}

In a source-free background,
\[
j_3^{\mathrm{D6/O6}}=0,
\]
and hence
\begin{equation}
dF_2'
=
\left(
\frac{2}{27}f^2
-
\frac{2}{5}m^2
\right)
e^\Phi\operatorname{Re}\Omega.
\label{eq:source-free-Bianchi}
\end{equation}

\subsubsection*{Identification with the geometric flux}

By \eqref{eq:SU3-structure-normalization}, the supersymmetric
NS--NS ansatz \eqref{eq:susy-H-flux} implies
\begin{equation}
\frac{h}{\sqrt{\lambda_1\lambda_2\lambda_3}}
=\frac{2m}{5}e^\Phi .
\label{eq:h-Romans-dilaton-relation}
\end{equation}
The factor relating $\eta_0$ to $\operatorname{Re}\Omega$ depends
on the metric scale. The relation above is a constraint on the
supersymmetric ansatz. The prescribed and evolving flux laws in the
geometric problem are defined independently of that constraint.

\subsection{Relation with gauge data}

In Scherk--Schwarz reductions, the gauging is encoded by a constant
twist matrix. On a generalized flag manifold, the analogous
homogeneous algebraic information is distributed among the
components of the reductive bracket.

This analogy motivates the schematic correspondence
\[
M
\longleftrightarrow
\mu,
\]
but it does not identify an arbitrary bracket-flow trajectory with a
renormalization-group flow of a physical gauge theory. Such an
interpretation requires a consistent truncation in which the
homogeneous metric and flux parameters are identified with
lower-dimensional scalar fields and gauging data.

Accordingly, in this work the bracket flow is treated as a geometric
evolution. Its possible interpretation as an evolution of effective
gauge data is left as a question for a specified compactification
model.

\subsection{Curvature Functionals and Effective Potentials}
\label{subsec:curvature-effective-potential}

The bracket-flow equation is driven by a geometric operator,
\[
P_{\mu,H}
=
\operatorname{Ric}_{\mu}^{\sharp}
-
\frac14\mathcal K_{\mu,H}^{\sharp}.
\]
This does not, by itself, imply that the flow is Hamiltonian or that
it is the gradient flow of a lower-dimensional scalar potential.

For homogeneous Ricci flow, suitable normalized versions of the
bracket flow admit variational interpretations involving the scalar
curvature or the norm of the moment map. In the presence of a
three-form flux, however, the relevant functional must also contain
the flux energy.

A natural geometric quantity is
\begin{equation}
\mathcal F(g,H)
=
\int_M
\left(
\operatorname{Scal}_g
-
\frac{1}{12}H_{abc}H^{abc}
\right)
dV_g.
\label{eq:geometric-flux-functional}
\end{equation}
For an invariant structure on a compact homogeneous space, its
integrand is constant, and hence
\[
\mathcal F(g,H)
=
\operatorname{Vol}(M,g)
\left(
\operatorname{Scal}_g
-
\frac{1}{12}H_{abc}H^{abc}
\right).
\]

The lower-dimensional scalar potential obtained from a
compactification generally depends not only on the homogeneous
bracket and the flux, but also on the volume modulus, the dilaton,
the choice of conformal frame, and the normalization of the
higher-dimensional action. We therefore write schematically
\begin{equation}
V_{\mathrm{eff}}
=
V_{\mathrm{eff}}
(\boldsymbol\lambda,H,\varphi),
\label{eq:effective-potential-schematic}
\end{equation}
where $\varphi$ denotes any additional scalar fields arising in the
dimensional reduction.

Accordingly, a relation of the form
\[
\dot\phi^i
=
-\Omega(\phi)\,
G^{ij}\partial_jV_{\mathrm{eff}}
\]
must be verified after the scalar-field metric $G_{ij}$ and the
effective potential have been derived. Such a gradient-flow
interpretation will not be assumed in the general discussion.

\subsection{Geometric evolution and moduli-space distance}
\label{subsec:moduli-space-distance}

The diagonal invariant metrics are parametrized by
\[
\boldsymbol\lambda
=
(\lambda_1,\ldots,\lambda_s),
\qquad
\lambda_i>0.
\]
It is convenient to introduce logarithmic variables
\[
u_i=\log\lambda_i.
\]
A homogeneous geometric flow then determines a curve
\[
t\longmapsto
\boldsymbol u(t)
\]
in the space of invariant metrics.

Let $G_{ij}(\boldsymbol u)$ denote the metric induced on the scalar
sector of the lower-dimensional theory. The length of a trajectory
is
\begin{equation}
\Delta(t_0,t_1)
=
\int_{t_0}^{t_1}
\sqrt{
G_{ij}(\boldsymbol u(t))
\dot u_i(t)\dot u_j(t)
}\,dt.
\label{eq:moduli-distance}
\end{equation}

For an infinite-distance claim, one must show that the scalar-field
distance from a fixed initial point diverges in the metric of the
specified reduction. The length \eqref{eq:moduli-distance} is an
upper bound for that distance. Bracket degeneration alone gives
neither quantity. If the distance $d_{\mathrm{scal}}$ diverges, the
Swampland Distance Conjecture motivates a separate test of the
asymptotic mass law
\[
m_{\mathrm{tower}}\sim e^{-\alpha d_{\mathrm{scal}}},\qquad \alpha>0.
\] In the full flag
example, this requires determining both the asymptotic behavior
of the metric parameters and the corresponding Kaluza--Klein mass
scales.

\section{Conclusion}
\label{sec:conclusion}

The exact invariant three-form on $\mathrm{SU}(3)/T^2$ gives a
quadratic flux operator that acts by the same scalar on all three
isotropy summands. In the prescribed case, it balances Ricci curvature
at one metric on each invariant Einstein ray. The isotropic point has
unstable index two and the three K\"ahler--Einstein points have index
one. These equilibria also determine stationary complete reductive
brackets.

With Hodge heat evolution, the isotropic solution instead becomes
extinct in finite time with vanishing relative flux intensity.
The transverse variational equation has absolute decay of order
$(T-t)^{3/5}$ and relative growth of order $(T-t)^{-2/5}$.
The nonlinear dynamics of nearby anisotropic trajectories remain open.
A lower-dimensional interpretation of these geometric flows requires
a specified reduction, including the scalar kinetic metric, effective
potential and Kaluza--Klein spectrum.

\appendix
\section{Ricci normalization on the full flag manifold}
\label{app:ricci-su3}

We use the decomposition and $Q$-orthonormal bases of Section~\ref{subsec:su3-reductive-bracket}.

Since the isotropy summands are pairwise inequivalent, the Ricci
endomorphism is diagonal:
\[
\operatorname{Ric}_{\boldsymbol\lambda}^{\sharp}
|_{\mathfrak m_i}
=
r_i(\boldsymbol\lambda)\,I_{\mathfrak m_i}.
\]

For $Q$-orthonormal bases
\[
\left\{
e_{\alpha}^{(i)}
\right\}_{\alpha=1}^{d_i}
\subset\mathfrak m_i,
\]
define
\begin{equation}
[ijk]
=
\sum_{\alpha=1}^{d_i}
\sum_{\beta=1}^{d_j}
\sum_{\gamma=1}^{d_k}
Q\left(
[e_{\alpha}^{(i)},e_{\beta}^{(j)}]_{\mathfrak m_k},
e_{\gamma}^{(k)}
\right)^2.
\label{eq:structural-constants-ijk}
\end{equation}

For $\mathrm{SU}(3)/T^2$, one has
\[
d_1=d_2=d_3=2.
\]
The quantities $[ijk]$ are symmetric under permutations of their
indices by the invariance of $Q$. The four brackets displayed with
\eqref{eq:closed-flux-normalization} give $[123]=4c_0^2=1/3$.
With the normalization $Q=-B$, the only nonzero constants
with three pairwise distinct indices are
\begin{equation}
[123]=[132]=[213]=[231]=[312]=[321]=\frac13.
\label{eq:structural-constants-su3}
\end{equation}

The homogeneous Ricci formula gives
\begin{equation}
r_i
=
\frac{1}{2\lambda_i}
+
\frac{[123]}{2d_i}
\left(
\frac{\lambda_i}{\lambda_j\lambda_k}
-
\frac{\lambda_j}{\lambda_i\lambda_k}
-
\frac{\lambda_k}{\lambda_i\lambda_j}
\right),
\label{eq:ricci-general-su3}
\end{equation}
where $\{i,j,k\}=\{1,2,3\}$.

Substituting
\[
[123]=\frac13,
\qquad
d_i=2,
\]
into \eqref{eq:ricci-general-su3}, we obtain
\begin{equation}
\begin{aligned}
r_1
&=
\frac{1}{2\lambda_1}
+
\frac{1}{12}
\left(
\frac{\lambda_1}{\lambda_2\lambda_3}
-
\frac{\lambda_2}{\lambda_1\lambda_3}
-
\frac{\lambda_3}{\lambda_1\lambda_2}
\right),\\[0.3em]
r_2
&=
\frac{1}{2\lambda_2}
+
\frac{1}{12}
\left(
\frac{\lambda_2}{\lambda_1\lambda_3}
-
\frac{\lambda_1}{\lambda_2\lambda_3}
-
\frac{\lambda_3}{\lambda_1\lambda_2}
\right),\\[0.3em]
r_3
&=
\frac{1}{2\lambda_3}
+
\frac{1}{12}
\left(
\frac{\lambda_3}{\lambda_1\lambda_2}
-
\frac{\lambda_1}{\lambda_2\lambda_3}
-
\frac{\lambda_2}{\lambda_1\lambda_3}
\right).
\end{aligned}
\label{eq:ricci-components-su3}
\end{equation}

These quantities are the eigenvalues of the Ricci endomorphism,
rather than the coefficients of the Ricci tensor relative to the
fixed background metric $Q$. More precisely,
\[
\operatorname{Ric}_{g_{\boldsymbol\lambda}}
|_{\mathfrak m_i}
=
r_i(\boldsymbol\lambda)\,
g_{\boldsymbol\lambda}|_{\mathfrak m_i}
=
r_i(\boldsymbol\lambda)\lambda_iQ|_{\mathfrak m_i}.
\]

As a consistency check, consider the normal homogeneous metric
\[
\lambda_1=\lambda_2=\lambda_3=\lambda.
\]
Then
\[
r_1=r_2=r_3
=
\frac{5}{12\lambda}.
\]
Therefore,
\[
\operatorname{Ric}^{\sharp}
=
\frac{5}{12\lambda}I_{\mathfrak m},
\]
so the normal metric is Einstein, as expected.

For the unnormalized homogeneous Ricci flow,
\[
\frac{\partial}{\partial t}g
=
-2\operatorname{Ric}_g,
\]
the metric parameters satisfy
\begin{equation}
\dot\lambda_i
=
-2r_i(\boldsymbol\lambda)\lambda_i.
\label{eq:ricci-flow-lambda-su3}
\end{equation}
Thus,
\begin{equation}
\begin{aligned}
\dot\lambda_1
&=
-1
-\frac16
\left(
\frac{\lambda_1^2}{\lambda_2\lambda_3}
-
\frac{\lambda_2}{\lambda_3}
-
\frac{\lambda_3}{\lambda_2}
\right),\\
\dot\lambda_2
&=
-1
-\frac16
\left(
\frac{\lambda_2^2}{\lambda_1\lambda_3}
-
\frac{\lambda_1}{\lambda_3}
-
\frac{\lambda_3}{\lambda_1}
\right),\\
\dot\lambda_3
&=
-1
-\frac16
\left(
\frac{\lambda_3^2}{\lambda_1\lambda_2}
-
\frac{\lambda_1}{\lambda_2}
-
\frac{\lambda_2}{\lambda_1}
\right).
\end{aligned}
\label{eq:ricci-flow-explicit-su3}
\end{equation}

\section{Isotropic system with a general structural coefficient}
\label{app:isotropic-extinction}

Consider the following auxiliary ODE with a free nonzero
coefficient $c_0$ in the flux equation and the metric equation fixed
as on the full flag. The auxiliary system is

\begin{equation}
\begin{aligned}
\dot\lambda
&=-\frac56+\frac{2h^2}{\lambda^2},\\
\dot h
&=-\frac{12c_0^2}{\lambda}h,
\end{aligned}
\label{eq:isotropic-system-extinction}
\end{equation}
where $\lambda>0$ and $c_0\neq0$.

Introduce
\[
z=\frac{h^2}{\lambda^2},
\qquad
a=\frac53-24c_0^2,
\]
and define a new time parameter $\tau$ by
\[
\frac{d\tau}{dt}=\frac1{\lambda(t)},
\qquad \tau(t_0)=0.
\]
Writing a prime for differentiation with
respect to $\tau$, the system becomes
\begin{equation}
\begin{aligned}
z'&=z(a-4z),\\
(\log\lambda)'&=-\frac56+2z,\\
h'&=-12c_0^2h.
\end{aligned}
\label{eq:rescaled-isotropic-system}
\end{equation}

\begin{proposition}
\label{prop:explicit-isotropic-solutions}
Let $\lambda_0>0$, $h_0\in\mathbb R$, and
$z_0=h_0^2/\lambda_0^2$.

If $a\neq0$, the solution is given by
\begin{equation}
z(\tau)=
\frac{az_0e^{a\tau}}
{a+4z_0(e^{a\tau}-1)},
\label{eq:explicit-z}
\end{equation}
\begin{equation}
\lambda(\tau)=
\lambda_0e^{-5\tau/6}
\left[
1+\frac{4z_0}{a}(e^{a\tau}-1)
\right]^{1/2},
\label{eq:explicit-lambda}
\end{equation}
and
\begin{equation}
h(\tau)=h_0e^{-12c_0^2\tau}.
\label{eq:explicit-h}
\end{equation}

If $a=0$, the corresponding expressions are
\begin{equation}
z(\tau)=\frac{z_0}{1+4z_0\tau},
\label{eq:critical-z}
\end{equation}
and
\begin{equation}
\lambda(\tau)=
\lambda_0e^{-5\tau/6}
(1+4z_0\tau)^{1/2}.
\label{eq:critical-lambda}
\end{equation}
These formulas are valid for all $\tau\geq0$,
with the zero-flux case interpreted as
$z\equiv0$.
\end{proposition}

\begin{proof}
The first equation in
\eqref{eq:rescaled-isotropic-system}
is the logistic equation
\[
z'=az-4z^2.
\]
Direct integration yields
\eqref{eq:explicit-z} when $a\neq0$
and \eqref{eq:critical-z} when $a=0$.

Next, integrating
\[
(\log\lambda)'=-\frac56+2z
\]
and substituting the expression for $z$
gives \eqref{eq:explicit-lambda} and
\eqref{eq:critical-lambda}, respectively.

Finally, integration of
$h'=-12c_0^2h$ gives
\eqref{eq:explicit-h}.
\end{proof}

\begin{theorem}
\label{thm:isotropic-finite-extinction}
Every solution of
\eqref{eq:isotropic-system-extinction}
with $\lambda_0>0$ and $c_0\neq0$
has a finite maximal forward time $T$,
and
\[
\lim_{t\uparrow T}\lambda(t)=0,
\qquad
\lim_{t\uparrow T}h(t)=0.
\]

More precisely:

\begin{enumerate}
\item If $a\leq0$, then
\[
\frac{h(t)^2}{\lambda(t)^2}\longrightarrow0,
\qquad
\lambda(t)\sim\frac56(T-t).
\]

\item If $a>0$ and $h_0\neq0$, then
\[
\frac{h(t)^2}{\lambda(t)^2}
\longrightarrow\frac a4
=\frac5{12}-6c_0^2,
\]
and
\[
\lambda(t)\sim12c_0^2(T-t).
\]
\end{enumerate}

If $h_0=0$, the solution instead satisfies
$\lambda(t)=\lambda_0-\frac56(t-t_0)$,
regardless of the sign of $a$.
\end{theorem}

\begin{proof}
Since
\[
\frac{dt}{d\tau}=\lambda(\tau),
\]
the forward existence time in the original
parameter is
\[
T=t_0+\int_0^\infty\lambda(\tau)\,d\tau.
\]

Suppose first that $a<0$. From the explicit
solution,
\[
z(\tau)\longrightarrow0,
\qquad
\lambda(\tau)=O(e^{-5\tau/6}).
\]
The integral defining $T$ therefore
converges.

If $a=0$, then
\[
z(\tau)\longrightarrow0,
\qquad
\lambda(\tau)
=O\left(e^{-5\tau/6}\sqrt{1+\tau}\right),
\]
and the integral again converges.

Finally, if $a>0$ and $z_0>0$, then
\[
z(\tau)\longrightarrow\frac a4.
\]
Moreover,
\[
\lambda(\tau)
\sim C e^{(-5/6+a/2)\tau}
=C e^{-12c_0^2\tau}
\]
for some $C>0$.
Since $c_0\neq0$, this is integrable.
The case $z_0=0$ is immediate.

In every case,
$\lambda(\tau)\to0$ and $h(\tau)\to0$,
and the original time approaches a finite
limit $T$.

For $a\leq0$, we have $z(t)\to0$.
Consequently,
\[
\dot\lambda(t)\longrightarrow-\frac56.
\]
Since $\lambda(T)=0$ in the limiting sense,
integration gives
\[
\lambda(t)\sim\frac56(T-t).
\]

For $a>0$ and $z_0>0$,
\[
\dot\lambda(t)
\longrightarrow
-\frac56+\frac a2
=-12c_0^2.
\]
Therefore,
\[
\lambda(t)\sim12c_0^2(T-t).
\]
This completes the proof.
\end{proof}

\begin{remark}
The parameter
\[
a=\frac53-24c_0^2
\]
determines the asymptotic behavior of the
relative flux intensity. For the full flag with $Q=-B$, equation
\eqref{eq:closed-flux-normalization} gives $a=-1/3$.

The extinction established here is
isotropic: all six tangent directions
shrink at the same rate. It should not
be interpreted as a selective dimensional
compactification.
\end{remark}

\section{Scherk--Schwarz twists and root data}

In a circle reduction, a Scherk--Schwarz twist is represented by
$U(y)=\exp(yM)$, with $M$ in the Lie algebra of the duality group
\cite{ScherkSchwarz1979,Hull2002,Samtleben2008}. In a homogeneous
flag, the relevant geometric coefficients are the components of
the reductive bracket. For a Cartan--Weyl basis, the root-sum
components are $[X_\alpha,X_\beta]=N_{\alpha,\beta}X_{\alpha+\beta}$
when $\alpha+\beta$ is a root. Opposite roots contribute to the
Cartan component. The scalar coefficients of an invariant three-form
have weight zero and obey the separate closure condition $d_\mu H=0$.

For example, on $\mathrm{SU}(4)/T^3$ the positive roots are
\[
\alpha_1,\ \alpha_2,\ \alpha_3,\ \alpha_1+\alpha_2,\
\alpha_2+\alpha_3,\ \alpha_1+\alpha_2+\alpha_3 .
\]
Relations such as
$\alpha_1+(\alpha_2+\alpha_3)=\alpha_1+\alpha_2+\alpha_3$
describe the corresponding root brackets. The resulting constants
are fixed by the Lie algebra and the chosen basis. Invariant forms
can be compared with generalized flux constructions
\cite{Hitchin2003,Gualtieri2011,SheltonTaylorWecht2005,
AldazabalMarquesNunez2006,Blumenhagen2015,GasparimGramaVarea2023}.
Identifying lower-dimensional gauging data requires a specified
reduction, such as the invariant coset reduction of
\cite{CassaniKashaniPoor2009}.


\begin{thebibliography}{99}

\bibitem{Caviezel2008}
C.~Caviezel, P.~Koerber, S.~K\"ors, D.~L\"ust, D.~Tsimpis and M.~Zagermann,
The effective theory of type IIA AdS$_4$ compactifications on nilmanifolds and cosets,
\emph{Classical and Quantum Gravity} \textbf{26}, 025014 (2009).
\href{https://doi.org/10.1088/0264-9381/26/2/025014}{doi:10.1088/0264-9381/26/2/025014}.

\bibitem{GasparimGramaVarea2023}
E.~Gasparim, L.~Grama and C.~Varea,
The $H$-flux on flag manifolds generated by infinitesimal $T$-duality,
\emph{Letters in Mathematical Physics} \textbf{116}, 72 (2026).
\href{https://doi.org/10.1007/s11005-026-02106-3}{doi:10.1007/s11005-026-02106-3}.

\bibitem{ScherkSchwarz1979}
J. Scherk and J. H. Schwarz, How to get masses from extra dimensions, Nucl. Phys. B, vol. 153, pp. 61--88, 1979.

\bibitem{Hull2002}
C. M. Hull, Gauged d=9 supergravities and scherk-schwarz reduction, Classical and Quantum Gravity, vol. 21, pp. 509--516, 2004.

\bibitem{Samtleben2008}
H. Samtleben, Lectures on gauged supergravity and flux compactifications, Class. Quant. Grav., vol. 25, p. 214002, 2008.

\bibitem{PediconiSbiti2022}
F.~Pediconi and S.~Sbiti,
Collapsed ancient solutions of the Ricci flow on compact homogeneous spaces,
\emph{Proceedings of the London Mathematical Society} \textbf{125}, 1130--1151 (2022).
\href{https://doi.org/10.1112/plms.12478}{doi:10.1112/plms.12478}.

\bibitem{BaumePalti2016}
F. Baume and E. Palti, Backreacted axion field ranges in string theory, JHEP, vol. 1608, p. 043, 2016.

\bibitem{CalderonInfante2020}
J.~Calder\'on-Infante, A.~M. Uranga and I.~Valenzuela,
The convex hull swampland distance conjecture and bounds on non-geodesics,
\emph{JHEP} \textbf{03} (2021), 299.

\bibitem{LustPaltiVafa2019}
D. Lust, E. Palti, and C. Vafa, Ads and the swampland, Physics Letters B, vol. 797, p. 134867, 2019.

\bibitem{Kehagias2020}
A. Kehagias, D. Lust, and S. Lust, Swampland, ricci flow and distance conjectures, Fortschritte der Physik, vol. 69, p. 2100096, 2021.

\bibitem{DemulderRaml2025}
S.~Demulder, D.~L\"ust and T.~Raml, Navigating string theory field space with geometric flows,
\emph{JHEP} \textbf{05} (2025), 030.

\bibitem{Lauret2013CurvatureFlowsHomogeneous}
J. Lauret, Ricci flow of homogeneous manifolds, Mathematische Zeitschrift, vol. 274, no. 1--2, pp. 373--403, 2013.

\bibitem{Lauret2011BracketFlow}
J. Lauret, The ricci flow for simply connected nilmanifolds, Communications in Analysis and Geometry, vol. 19, no. 5, pp. 831--854, 2011.

\bibitem{GarciaFernandez2013}
M. García-Fernández, Torsion-free generalized connections and heterotic supergravity, Communications in Mathematical Physics, vol. 332, no. 1, pp. 89--115, 2014.

\bibitem{StreetsTian2010}
J. Streets and G. Tian, A parabolic flow of pluriclosed metrics, International Mathematics Research Notices, vol. 2010, no. 16, pp. 3101--3133, 2010.

\bibitem{FusiLafuenteStanfield2024}
E.~Fusi, R.~A.~Lafuente and J.~Stanfield,
The homogeneous generalized Ricci flow,
\href{https://arxiv.org/abs/2404.15749}{arXiv:2404.15749} (2024).



\bibitem{GarciaFernandezStreets2021}
M.~Garc\'ia-Fern\'andez and J.~Streets,
\emph{Generalized Ricci Flow}, University Lecture Series, vol.~76,
American Mathematical Society, 2021.

\bibitem{Hitchin2003}
N. Hitchin, Generalized calabi-yau manifolds, Q. J. Math., vol. 54, no. 3, pp. 281--308, 2003.

\bibitem{Gualtieri2011}
M. Gualtieri, Generalized complex geometry, Ann. of Math. (2), vol. 174, no. 1, pp. 75--123, 2011.

\bibitem{SheltonTaylorWecht2005}
J. Shelton, W. Taylor, and B. Wecht, Non-geometric flux compactifications, JHEP, vol. 10, p. 085, 2005.

\bibitem{AldazabalMarquesNunez2006}
G. Aldazabal, D. Marques, and C. Núñez, Non-geometric fluxes and gauged supergravities, JHEP, vol. 11, p. 052, 2006.

\bibitem{Blumenhagen2015}
R. Blumenhagen, A. Deser, and E. Plauschinn, Non-geometric fluxes and their algebra, Fortsch. Phys., vol. 63, pp. 619--643, 2015.

\bibitem{CassaniKashaniPoor2009}
D.~Cassani and A.-K.~Kashani-Poor,
Exploiting $\mathcal N=2$ in consistent coset reductions of type IIA,
\emph{Nuclear Physics B} \textbf{817}, 25--57 (2009).

\end{thebibliography}
\end{document}